\documentclass[11pt]{article}
\usepackage{fullpage}
\usepackage{graphics}
\usepackage{amsfonts}
\usepackage{times}
\begin{document}
\title{Queues with heterogeneous servers and uninformed customers: who works the most?}
\author{Fabricio Bandeira Cabral\\
Universidade Federal de Pelotas, UNIPAMPA\\
 fabriciocabral2000@yahoo.com }

\date{ }
\maketitle
\newtheorem{Theo}{Theorem}[section]
\newtheorem{Lemma}{Lemma}[section]
\newtheorem{Fact}{Lemma}[section]
\newtheorem{Cor}{Corollary}[section]
\newtheorem{Prop}{Property}[section]
\newtheorem{Def}{Definition}[section]
\newcommand{\bx}{\rightline{$\Box$}}
\newcommand{\Proof}{{\em Proof.\ }}
\newcommand{\definedas}{\stackrel{\Delta}{=}}
\newcommand{\And}{ \, \, \, \mbox{ and }}
\newcommand{\with}{ \, \, \, \mbox{ with }}

\begin{abstract}
 In this paper, we consider systems that can
be modelled by $M \mid M \mid n$ queues with  heterogeneous
 servers and non informed customers. Considering any two servers: we show that the probability
 that the fastest server is busy is smaller than the probability
 that the  slowest server is busy. Moreover, we show that the effective rate of service done by the fastest server
is larger than effective rate of service done by the slowest server.
\end{abstract}
~

\section{Introduction}

In this work, we consider queues with heterogeneous server, we focus
on the uninformed customers case~\cite{Rubi,Cabral}. A practical
motivation for this is that there are plenty of systems with servers
which operate in parallel and are heterogeneous in their
capabilities, but without the customer being aware of their service
rates differences. For instance: tellers in a bank, cashiers in a
supermarket, agents in airport check-in, and may others.

 Thus, in order to determine who works the most
for uninformed customers, we model the system by a $M \mid M \mid n$
queue~\cite{Kleinrock} with heterogeneous servers and uninformed
customers. The relevant probability distribution obtained through
the use of the balance equations was derived in \cite{Cabral}.

Considering any two servers $l$ and $m,$ and using the probability
distribution: we compute the difference between the probability that
server $l$ is busy and the probability that server $m$ is busy:
\begin{itemize}

\item we compute the difference between the probability that
server $l$ is busy and the probability that server $m$ is busy
$\mathcal{P}_{l}^n - \mathcal{P}_{m}^n;$

\item we compute the difference between the effective rates of  service done by
servers  $l$ and $m$ is
$\mu_l\mathcal{P}_{l}^n-\mu_m\mathcal{P}_{l}^n.$

\end{itemize}
For  a given $0<\lambda<\mu^n,$ and $\mu_l > \mu_m,$ the first
difference is negative and second difference is positive, which
leads to the desired result.

The third expression given in the theorem only combines these two
results in an upper bound -lower bound form.

\section{ $M \mid M \mid n$ queues with heterogeneous servers and uninformed
customers \label{MMN}} The system for uninformed customers under
Markovian assumptions can be modelled by a Markov process.
The states of this Markov process are associated to $(n+1)$-tuples
$(x_1,x_2,\ldots,x_n,x_f),$ where $x_i,$  for $i=1$ to $n,$ is a
boolean variable, which is  equal to one iff server $i$ is busy,
and $x_f$ denotes the length of the queue of the system. We shall
index the states by
\[i= x_1 2^{n-1}+x_2 2^{n-2}+\ldots+x_n 2^0+x_f,\]
where
\[ x_f> 0 \Rightarrow x_1=x_2=\ldots=x_n=1. \]

\subsection{Notation}
 In order to study the steady state behavior
of the Markov process, we utilize the  notation defined in
table~\ref{tabnot}.

\begin{table}
\caption{Table of notations}
 \label{tabnot}
\begin{tabular}{|l|l|}
  \hline
  $U_k^n$ & set of states with $k$ users for the system  with $n$
  servers \\
  $S_i^n$ & set of busy servers when the system with $n$
  servers  is in state $i$\\
  neighbor states & two states $i$ and $j$  that differ by exactly \\
           & one component of the $(n+1)$-tuple for the system with $n$
  servers \\
  $g^n(i,j)$ &  for two states $i$ and $j$  that differ by only \\
           & one component of the $(n+1)$-tuple (neighbor
           states),\\
           & $g^n(i,j)$ denotes the index of that component, \\
                     &for the system with $n$
  servers  \\
  $D^n$ & set of ordered pairs of neighbor states for the system with $n$
  servers \\
  $D_i^n$ & the set of neighbor states of state $i$ for the system with $n$
  servers \\
  $\mu_l$ & average service rate for server $l$ \\
$\mu_1 \geq \mu_2 \geq \ldots \geq \mu_n$&\\
  $\lambda$ & arrival rate \\
  $q(i,j)$ & transition rate between states $i$ and $j$\\
  $N^n$ & average number of customers in the system with $n$
  servers \\
  $T^n$ & average sojourn time of a customer in the system with $n$
  servers \\
  $\mu^n$  & $\sum_{i=1}^n  \mu_i$ \\
  $p_i^n$ & steady state probability $p_i^n$
             that the system with $n$
  servers is in state $i$\\
  $P_k^n$   &        probability that the number of users in the system with $n$
  servers is $k$\\
  $\mathcal{P}_l^n$   &        probability that server $l$ is busy in the system with $n$
  servers is $k$\\
  $\mathcal{P}_{l,m}^n$   &        probability that server $l$ is busy and server $m$ is idle in the system with $n$
  servers \\
$\mathcal{P}^n$   &        probability that all servers $l$ are busy
in the system with $n$
  servers \\
  $S^n=\{\mu_1,\mu_2,\ldots,\mu_n\}$&\\
  $S^{-n}=\{\mu_1^{-1},\mu_2^{-1},\ldots,\mu_n^{-1}\}$&\\
 $P(\mathcal{S},k)$ &  denotes the sum of all distinct products of $k$ elements
 of\\
 & the set $\mathcal{S}$ \\
&\\
\hline
\end{tabular}
\end{table}

\subsection{Balance Equations}
\subsubsection{Preliminary Considerations}
 As in \cite{Cabral}, let us consider a state $i$ for a system with
$n$ servers. At state $i,$ the system has $\mid S_i^n \mid$ busy
servers.
 The states that are neighbors of  state $i$ are given by the set
$D_i^n.$ There are transitions to state $i,$ from states in which
the system has $\mid S_i^n \mid - 1$ busy servers and from states
in which the system has $\mid S_i^n \mid + 1$ busy servers.

The set of states that are neighbor of  state $i$  and in which
the system has $\mid S_i^n \mid + 1$ busy servers is given by
$U^n_{\mid S_i^n \mid +1} \cap D_i^n.$ The transitions from a
state $j$ in this set to state $i$ happen when server $g^n(i,j)$
finishes serving its customer.

The set of states that are neighbor of  state $i$  and in which the
system has $\mid S_i^n \mid - 1$ busy servers is given by $U^n_{\mid
S_i^n \mid - 1} \cap D_i^n.$ The transitions from a state $j$ in
this set to state $i$ happen when a customer arrives at the system
in state $j.$ This means that one of the $n-(\mid S_i^n \mid -1)$
idle servers will become busy. There is no preference among these
servers, which implies that the transition rate is given by the
arrival rate divided by the number of idle servers.

\subsubsection{Balance Equations} The steady state probability
$p_i^n$, probability that the system is in state $i,$ is given by
the solution to the following linear system:

\[
(\lambda + \sum_{k \in S_i^n} \mu_k)p_i^n -\sum_{j \in U^n_{\mid
S_i^n \mid -1} \cap D_i^n}\frac{\lambda}{n-(\mid S_i^n \mid
-1)}p_j^n -\sum_{j \in U^n_{\mid S_i^n \mid +1} \cap D_i^n}
\mu_{g^n(i,j)} p_j^n = 0,\]
\[0 \leq i \leq 2^n -1;\]

\[\lambda p_{i-1}^n=\mu^n p_{i}^n, \, i \geq 2^n;\]
\[\sum_{i=0}^\infty  p_i^n=1.\]

\subsubsection{Main Results}
\begin{Theo}
Let us assume that the Markov process with which we are dealing is
time homogeneous, irreducible, and remains in each state for a
positive length of time and is incapable of passing through an
infinite number of states in a finite time~\cite{Kelly}. Then, for a
given $0<\lambda<\mu^n,$ and $\mu_l > \mu_m,$
we have that
\begin{itemize}

\item $\mathcal{P}_{l}^n < \mathcal{P}_{m}^n;$

\item $\mu_l\mathcal{P}_{l}^n > \mu_m \mathcal{P}_{m}^n;$

\item $ \frac{\mu_m}{\mu_l}\mathcal{P}_{m}^n < \mathcal{P}_{l}^n < \mathcal{P}_{m}^n.$
\end{itemize}

 \label{F1}
\end{Theo}

\Proof

 From~\cite{Cabral}, we have that
\[
p_i^n=\frac{(n-\mid S_i^n \mid)!}{n!} \frac{\lambda^{\mid S_i^n
\mid}}{\prod_{j \in S_i^n} \mu_j} p_0^n,  \,\,\, 0 \leq i \leq 2^n -
1.
\]

For a state  $i, \,\,\, 0 \leq i \leq 2^n - 1,$ such that server $l$
is busy, server $m$ is idle,and there are other $k-1$ servers busy,
we have that

\begin{eqnarray*}
p_i^n &=& \frac{\lambda}{\mu_l} \frac{(n-k)!}{n!}
\frac{\lambda^{k-1}}{\prod_{j \in S_i^n -
\{\mu_l,\mu_m\} } \mu_j} p_0^n\\
\end{eqnarray*}

The  probability that server $l$ is busy and server $m$ is idle in
the system with $n$ servers is then given by
\[\mathcal{P}_{l,m}^n
=\frac{\lambda}{\mu_l} \sum_{k=1}^{n-1}
\frac{(n-k)!}{n!}\lambda^{k-1} P(S^{-n} -
\{\mu_l^{-1},\mu_m^{-1}\},k-1) p_0^n.
\]

 The difference between the
probability that server $l$ is busy and the probability that server
$m$ is busy is then given by
\[\mathcal{P}_{l,m}^n - \mathcal{P}_{m,l}^n
=(\frac{\lambda}{\mu_l}- \frac{\lambda}{\mu_m})\sum_{k=1}^{n-1}
\frac{(n-k)!}{n!}\lambda^{k-1} P(S^{-n} -
\{\mu_l^{-1},\mu_m^{-1}\},k-1) p_0^n. ,
\]

The probability that server $l$ is busy is given by
\[\mathcal{P}_{l}^n
=\frac{\lambda}{\mu_l} \sum_{k=1}^{n-1}
\frac{(n-k)!}{n!}\lambda^{k-1} P(S^{-n} - \{\mu_l^{-1}\},k-1) p_0^n
+ \mathcal{P}^n.
\]

The effective rate of of service done by server  $l$ is given by
$\mu_l \mathcal{P}_{l}^n.$

Thus, the difference between the effective rates of  service done by
servers  $l$ and $m$ is given by
\[\mu_l\mathcal{P}_{l}^n-\mu_m\mathcal{P}_{l}^n
=\frac{\lambda}{\mu_l} \sum_{k=1}^{n-1}
\frac{(n-k)!}{n!}\lambda^{k-1} (P(S^{-n} - \{\mu_l^{-1}\},k-1)-
P(S^{-n} - \{\mu_m^{-1}\},k-1)) p_0^n + (\mu_l-\mu_m)\mathcal{P}^n.
\]

Noticing that
\[P(S^{-n}- \{\mu_l^{-1}\},k-1)=\frac{1}{\mu_m} P(S^{-n} -
\{\mu_l^{-1},\mu_m^{-1}\},k-2) +  P(S^{-n} -
\{\mu_l^{-1},\mu_m^{-1}\},k-1).
\]
We have that
\[\mu_l\mathcal{P}_{l}^n-\mu_m\mathcal{P}_{l}^n= (\frac{1}{\mu_m}-\frac{1}{\mu_l})(\sum_{k=1}^{n-1}
\frac{(n-k)!}{n!}\lambda^{k} P(S^{-n} -
\{\mu_l^{-1},\mu_m^{-1}\},k-2)) + (\mu_l - \mu_m)\mathcal{P}.
\]

 \bx

\section{Applications \label{Appl}}
There are several systems with servers  working in parallel.
Moreover in many of these systems, the servers are heterogeneous in
their capabilities, without the customers being aware of who is fast
and who is slow.As examples of systems of this kind, we could cite:
tellers in a bank, cashiers in a supermarket, agents in an airport
check-in, etc~\cite{Cabral}.

\section{Concluding remarks \label{Conc}}
 In order to answer the question: who works the most?
 We cosidered systems that could
be modelled by $M \mid M \mid n$ queues with  heterogeneous
 servers and non informed customers. Through the use of the
 probability distribution previously derived, we were able to show
 that the fastest server works less in the sense that the
 probability that the fastest server is busy is smaller than the probability
 that the  slowest server is busy. On the
 other hand, we showed that effective rate of service of the fastest
 server is larger than that of the slowest server. These results
 were also combined in an upper bound - lower bound form.

\end{document}